\documentclass[12pt,reqno]{amsart}

\usepackage[T1]{fontenc}
\usepackage[utf8]{inputenc}
\usepackage{lmodern}
\usepackage{microtype}
\usepackage{amsmath,amssymb,amsthm,amsfonts,mathtools}
\usepackage{mathrsfs}
\usepackage{enumitem}
\usepackage{txfonts}
\usepackage{xcolor}
\usepackage[colorlinks,final,backref=page,hyperindex]{hyperref}
\usepackage{cleveref}

\setlist{itemsep=3pt,topsep=5pt}

\newtheorem{theorem}{Theorem}[section]
\newtheorem{proposition}[theorem]{Proposition}
\newtheorem{lemma}[theorem]{Lemma}
\newtheorem{corollary}[theorem]{Corollary}
\theoremstyle{definition}
\newtheorem{definition}[theorem]{Definition}
\newtheorem{example}[theorem]{Example}
\theoremstyle{remark}
\newtheorem{remark}[theorem]{Remark}

\newcommand{\bfk}{\mathbf{k}}
\newcommand{\Ga}{\mathbb G_{\mathrm a}}
\newcommand{\Gm}{\mathbb G_{\mathrm m}}
\newcommand{\GL}{\mathrm{GL}}
\newcommand{\cO}{\mathcal O}
\newcommand{\Dist}{\operatorname{Dist}}
\newcommand{\Lie}{\operatorname{Lie}}
\newcommand{\End}{\operatorname{End}}

\newcommand{\Hom}{\operatorname{Hom}}
\newcommand{\Spec}{\operatorname{Spec}}
\newcommand{\Ad}{\operatorname{Ad}}

\newcommand{\id}{\mathrm{id}}
\newcommand{\RB}{\operatorname{RB}}
\newcommand{\DEnd}{\operatorname{DEnd}}

\newcommand{\pr}{\operatorname{pr}}

\newcommand{\bB}{\mathcal B}

\newcommand{\cM}{\mathcal M}
\newcommand{\widehatG}{\widehat G_e}

\title[Rota--Baxter operators on algebraic groups]{Rota--Baxter operators on algebraic groups: algebraization, cocycles, and moduli}

\author{Chengming Bai}
\address{Chern Institute of Mathematics and LPMC, Nankai University, Tianjin, 300071, China}
\email{baicm@nankai.edu.cn}
\author{Li Guo}
\address{Department of Mathematics and Computer Science, Rutgers University, Newark, NJ, 07102, USA}
\email{liguo@rutgers.edu}
\author{Chuijia Wang}
\address{Department of Mathematics, The University of Hong Kong, Pokfulam, Hong Kong, China}
\email{wangchuijia@u.nus.edu}


\begin{document}

\begin{abstract}
In this article, we study two fundamental problems for Rota--Baxter operators on algebraic groups: the algebraization problem and the moduli problem. To do this, we introduce several new concepts in this subject. We characterize when a Rota--Baxter operator on a Lie algebra comes from that on a connected algebraic group based on the classical Hopf-theoretic reconstruction of algebraic groups. We introduce gauge actions on the set of Rota--Baxter operators on a given algebraic group, the restriction of which to the invertible locus intertwines the coboundary action on the set of invertible 1-cocycles. Such actions have simple descriptions in terms of graph realizations of Rota--Baxter operators. Together with a description in terms of difference-endomorphisms, we describe the moduli of invertible Rota--Baxter operators. For connected reductive groups, based on techniques developed by Brion, we show that the corresponding good moduli space for the invertible locus is discrete, which has various applications to post-algebraic groups and factorization methods. 
\end{abstract}

\maketitle

\vspace{-.8cm}

\tableofcontents

\vspace{-.8cm}

\section{Introduction}

Rota--Baxter operators on groups connect Lie-theoretic factorization and integrable systems with post-groups and solutions of the Yang--Baxter equation. In this paper, we study these operators on \emph{affine algebraic groups}, where the operator is a regular morphism of schemes and all induced structures are formulated in the algebraic-geometric context. The explicit identity for Rota--Baxter operators is already known for abstract and Lie groups. Furthermore, the corresponding formulation for affine group schemes also appears in the literature. However, in this article, our purpose is different and has two aspects. Namely, for a \emph{fixed} algebraic group, we want to determine when an infinitesimal Rota--Baxter operator comes from a global one and to describe the classification of invertible operators under the gauge equivalences in algebraic families. To do this, we first recall the Lie-theoretic background and then study its variants in the setting of affine algebraic groups.

\subsection{From modified classical Yang-Baxter equations to Rota--Baxter operators on Lie groups}

In the Lie-theoretic setting, the classical approach to Hamiltonian Lax equations in terms of $r$-matrix leads to factorization methods for
integrable systems. In the operator form, a weight-one Rota--Baxter
operator on a Lie algebra $\mathfrak g$ is a linear map $R$ satisfying
\[[R(x),R(y)]
 =R\bigl([R(x),y]+[x,R(y)]+[x,y]\bigr),\quad\forall x,y\in \mathfrak g.\]
Equivalently, the operator $A=2R+\id_{\mathfrak g}$ satisfies the normalized \emph{modified classical Yang--Baxter equation}
\[[A(x),A(y)]-A\bigl([A(x),y]+[x,A(y)]\bigr)=-[x,y],\quad\forall x,y\in \mathfrak g.\]
By Semenov-Tian-Shansky's factorization method \cite{STS}, these operators are closely related to the integration of appropriate Lax flows. Such a connection is already mentioned in the formulation of Guo--Lang--Sheng for Lie groups\cite{GLS}.

More precisely, let $G$ be a Lie group or more generally an abstract group. Guo, Lang and Sheng \cite{GLS} introduced the concept of Rota--Baxter operators on $G$, by the identity
\[
 B(g)B(h)=B\bigl(gB(g)hB(g)^{-1}\bigr),\quad \forall g,h\in G.
\]
For a Lie group $G$ with a smooth Rota--Baxter operator $B$, one can define another Lie-group structure $\circ_B$ on $G$ by
\[g\circ_Bh:=gB(g)hB(g)^{-1},\quad\forall g,h\in G.\]
Meanwhile, the maps
\[B_-(g)=B(g)\text{ and }B_+(g)=gB(g)\]
define homomorphisms of Lie groups from $(G,\circ_B)$ to $G$, which lead to the factorization
$g=B_+(g)B_-(g)^{-1}$. By differentiating at the identity, one gets the weight-one Rota--Baxter identity on the Lie algebra. Thus, the construction for Lie groups can be regarded as a global geometric counterpart of the operator
factorization.

By the work in\cite{BG,BGST}, these operators are further related to \emph{set-theoretic solutions to the Yang--Baxter equation} . For instance, one can define a post-group structure $\triangleright_B$ on $G$ by 
\[g\triangleright_Bh:=B(g)hB(g)^{-1},\quad\forall g,h\in G.\]
Moreover, these two group laws lead to a bijective Yang--Baxter map in the braid form. 

\subsection{From Lie groups to algebraic groups: algebraic Yang--Baxter maps}

The concept of integrability and Yang--Baxter formulations also have several important algebraic-geometric counterparts. The profound work of Hitchin on algebraically completely integrable systems arising from moduli spaces
of bundles illustrates the important role of abelian varieties in the algebraic integrability problem \cite{Hitchin}. On the other hand, Kogan and Zelevinsky studied integrable systems on certain symplectic leaves of complex semisimple Poisson--Lie groups \cite{KoganZelevinsky}.  Moreover, geometric crystals and matrix factorizations yield \emph{birational} Yang--Baxter maps on algebraic varieties \cite{EtingofGeom,Veselov}. These developments  motivate our study of Rota--Baxter operators in the context of algebraic groups.

Let $G$ be an affine algebraic group with the Lie algebra $\mathfrak g$. For us, a \emph{Rota--Baxter operator}
on $G$ means a \emph{morphism of schemes} $B:G\to G$ satisfying the algebraic Rota--Baxter identity, which is regarded as an identity of morphisms from $G\times G$ to $G$. A necessary condition is that its functor of points satisfy the corresponding Rota--Baxter identity for abstract groups. The pair $(G,B)$ is called a \emph{Rota--Baxter algebraic group}.

We also define an analogue of post groups in algebraic geometry. To distinguish the post-group structures introduced here from the post-group structures on groups and Lie groups, we use the term \emph{post-algebraic group} for a post-group whose operations are morphisms of affine schemes (see \Cref{sec:applications}). 

For a Rota--Baxter operator $B$ on an affine algebraic group, the resulting factorization maps $B_\pm$ are morphisms of algebraic groups from $(G,\circ_B)$ to $G$. Moreover, the induced post-group structure and the braid-form Yang--Baxter maps are regular morphisms. In particular, the obtained solution of the Yang--Baxter equation is realized by a \emph{schematic automorphism} $G\times G\xrightarrow{\sim}G\times G$, not merely a set bijection.

Note that the notion of \emph{affine Rota--Baxter group schemes} have been introduced \cite{MZL} and their relation with co-Rota--Baxter Hopf algebras and affine skew braces have been studied. The passage to the corresponding regular Yang--Baxter constructions is also stressed there. Our work puts additional emphasis on the \emph{obstruction to extending infinitesimal operators to global ones} and on the \emph{variation of operators in algebraic families}. We call the above two fundamental problems the ``algebraization'' problem and the ``moduli'' problem.

By \emph{algebraization}, we mean the passage from Lie-theoretic data to a morphism of the corresponding connected algebraic group. More precisely, assume $G$ is connected. For a weight-one Rota--Baxter operator $R$ on $\mathfrak{g}$, we ask whether it is the differential at the identity of a \emph{globally regular} operator $B_R:G\to G$. Unlike in the local analytic and formal setting, a general algebraic group does not necessarily have a globally regular exponential map. Therefore, formal integration does not automatically provide a regular self-map. Indeed, the converse problem can be resolved in the following way. The first step is to construct a formal map at the identity with the help of the universal enveloping algebra, which is a purely algebraic construction. The second step is to capture the condition under which this formal map comes from a regular map on the chosen algebraic group via the injective algebra map $\iota_G:\cO(G)\longrightarrow U(\mathfrak g)^*$. It is direct to see that a necessary condition for algebraization is that the transpose map $\bB_R^\vee$ preserves $\iota_G(\cO(G))$, which is also sufficient when $G$ is connected (see \Cref{thm:algebraization}). The obstructions of the preservation condition are genuinely global. In particular, for tori, they can be characterized in terms of corresponding cocharacter lattices.

By \emph{moduli}, we mean the study of isomorphism classes of Rota--Baxter operators on a fixed connected algebraic group $G$. More precisely, we introduce certain gauge actions on the set of Rota--Baxter operators on $G$ and study the corresponding quotient stack and the good moduli space of the invertible locus  when $G$ is connected and reductive. In particular, these gauge actions correspond to the natural twisted actions determined by coboundaries under the dictionary between invertible Rota--Baxter operators and scheme-invertible nonabelian 1-cocycles. Using the difference characterization, the moduli problem can be reduced to the study of a locus of the endomorphism scheme determined by the invertible condition. We then apply Brion's rigidity result for endomorphism schemes to this locus to obtain our theorem for moduli.

\subsection{Main results: algebraization and moduli}

For the algebraization problem, we first adopt Goncharov's extension \cite{Goncharov} of $R$ to a Rota--Baxter operator $\bB_R$ on $U(\mathfrak g)$. For connected $G$, there is an embedding $\iota_G:\cO(G)\rightarrow U(\mathfrak g)^*$ given by the pairing between distributions and regular functions. We prove that the formal map determined by $R$ extends to a regular map $G\to G$ if and only if $\bB_R^\vee$ preserves this distinguished copy of $\iota_G(\cO(G))$. For split tori, this is exactly the condition that the induced linear map preserves the cocharacter lattice. We can see that the criterion detects an obstruction invisible in the Lie algebra alone.

For the classification problem, we first obtain a group theoretic characterization of Rota--Baxter operators. More precisely, we show that graphs of Rota--Baxter operators on $G$ are exact complements in $G\rtimes_{\Ad}G$. From this characterization, we construct a functorial gauge action and obtain an explicit admissibility condition for its infinitesimal counterpart. By taking the functional inversion, we can identify invertible operators with scheme-invertible conjugation 1-cocycles, so that the gauge action corresponds naturally to the twisted action defined by coboundaries. Moreover, these 1-cocycles in turn correspond to group endomorphisms $\theta$, for which $\delta_\theta(g)=\theta(g)g^{-1}$ is an isomorphism of schemes.

When $G$ is connected and reductive, Brion's rigidity theorem shows that there is an identification between connected components of $\underline{\Hom}_{\mathrm{gp}}(G,G)$ and corresponding conjugation orbits \cite{Brion}. We show that the difference-isomorphism locus 
$D_G$ corresponding to the invertible locus of Rota--Baxter operators is a
union of scheme-theoretic conjugation orbits. Our results also work for nonreduced families. For representatives $\theta_i$
of these components and $C_i=Z_G(\theta_i(G))$, we obtain a decomposition
\[
 [D_G/G]\simeq\coprod_i BC_i.
\]
Since the centralizers $C_i$ are linearly reductive, the quotient has a discrete good moduli space. This further implies that connected families of Rota--Baxter operators on a connected reductive algebraic group are locally constant modulo gauge actions in the fpqc topology. This rigidity property of moduli is one of our main new results. As byproducts, we also show that a connected family of associated strictly inner lifted post-algebraic groups and the corresponding factorizations must be constant modulo gauge actions.

The paper is organized as follows. In \Cref{sec:groups}, we introduce the definition of Rota--Baxter operators on algebraic groups,  characterize them in terms of the graph description, and study their differentiation. In \Cref{sec:hopf}, we study Rota--Baxter operators from the perspective of the universal enveloping algebra and formal duality. Moreover, we characterize the condition for algebraization. A dictionary between invertible Rota--Baxter operators and invertible $1$-cocycles is established in \Cref{sec:cocycles}, and is used to construct the quotient stack and its good moduli space in \Cref{sec:moduli}. In \Cref{sec:examples}, we give examples for split tori and connected simple algebraic groups. We also treat the nonreductive example $G=\Ga^r$. In \Cref{sec:applications}, we give applications of the moduli formulation to various notions of independent interest, including post-algebraic groups, post-Lie algebras, Yang--Baxter maps, and factorizations.

\section{Rota--Baxter operators on algebraic groups}\label{sec:groups}

We begin with the scheme-theoretic definition of a Rota--Baxter operator on a fixed algebraic group. We present a characterization of Rota--Baxter operators in terms of their graph description, apply this characterization to define a gauge action. Taking differentiation of the operator recovers a Rota--Baxter operator of weight one on the corresponding Lie algebra.

\subsection{Conventions and basic notions}

Throughout the paper, $\bfk$ is an algebraically closed field of characteristic zero. An algebraic group $G$ is a smooth affine group scheme of finite type over $\bfk$. All products and morphisms are over $\bfk$. We write $e$ for the identity
element and $\Ad_g$ for the inner automorphism of $G$ defined by
\[
\Ad_g(h)=ghg^{-1},\quad \forall g,h\in G.
\]
We will also use a similar notation $\Ad(g)$ for the adjoint action of $g\in G$ on $\mathfrak g$.

An invertible algebraic map means an isomorphism of schemes. A morphism between algebraic groups is \emph{pointed} if it carries the identity to the identity. Note that a pointed morphism need not be a group homomorphism.

\begin{definition}\label{def:group-rb}
Let $G$ be an algebraic group. A \emph{Rota--Baxter operator} on $G$ is a morphism of $\bfk$-schemes $B:G\to G$ such that
\begin{equation}\label{eq:group-rb}
  B(g)B(h)=B\bigl(g\Ad_{B(g)}(h)\bigr)
\end{equation}
as an identity of morphisms $G\times G\to G$.
\end{definition}

We also call $(G,B)$ a \emph{Rota--Baxter algebraic group}. The identity for abstract groups originates from \cite{GLS}, while the formulation for affine group schemes has been treated in \cite[Definition 3.4]{MZL}. We will restrict ourselves to the setting of smooth affine algebraic groups of finite type, which is suitable for differentiation, the algebraization theory, and the moduli problem. If $G$ is commutative, \eqref{eq:group-rb} means precisely that $B$ is a group endomorphism. For every $G$, the inversion map $g\mapsto g^{-1}$ and the constant map $g\mapsto e$ are Rota--Baxter operators on $G$. The relation with the concept of a Rota--Baxter operator on the corresponding universal enveloping algebra is established in \Cref{sec:hopf}.

The definition also works in families. For a $\bfk$-scheme $S$, let $\RB_G(S)$ be
the set of $S$-morphisms $G_S\to G_S$ satisfying \eqref{eq:group-rb} over $S$. This defines a functor on the category of $\bfk$-schemes.

\subsection{Graphs and exact complements}

Put
\[
  E=G\rtimes_{\Ad}G,\]
with the multiplication given by
\[ (x,a)(y,b)=\bigl(x\Ad_a(y),ab\bigr)\quad \forall x,y,a,b\in G.
\]
Denote by
\[
  K=\{e\}\times G,
  \quad N=G\times\{e\}.
\]
The subgroup $N$ is normal, while the subgroup $K$ is not necessarily normal. A closed subgroup scheme $H\subset E$ is said to be an \emph{exact complement} to
$K$ if the multiplication map $H\times K\to E$ is an isomorphism of $\bfk$-schemes.

\begin{proposition}\label{prop:graph}
For a morphism $B:G\to G$ of $\bfk$-schemes, the following conditions are equivalent.
\begin{enumerate}[label=\textup{(\roman*)}]
\item The map $B$ is a Rota--Baxter operator.
\item The graph 
\[
  \Gamma_B=\{(g,B(g)):g\in G\}
\]
of $B$ is a closed subgroup scheme of $E$.
\item The graph $\Gamma_B$ is an exact complement to $K$, that is, the multiplication map
\[
  \Gamma_B\times K\longrightarrow E
\]
is an isomorphism of $\bfk$-schemes.
\end{enumerate}
Every exact complement to $K$ arises from a unique Rota--Baxter operator.
\end{proposition}

\begin{proof}
For any two elements $(g,B(g))$ and $(h,B(h))$ in $\Gamma_B$, their multiplication is given by
\[
 (g,B(g))(h,B(h))
 =\bigl(g\Ad_{B(g)}(h),B(g)B(h)\bigr).
\]
Therefore, \eqref{eq:group-rb} is equivalent to the condition that $\Gamma_B$ is closed under multiplication. The same identity shows $B(e)=e$. The inverse of $g$ is given by
$h=\Ad_{B(g)^{-1}}(g^{-1})$. Hence, the graph is a subgroup scheme of $E$. It is closed because $G$ is separated.

The factorization and its inverse are explicitly given by
\[
 (g,a)=(g,B(g))\bigl(e,B(g)^{-1}a\bigr),
\]
\[
 (g,a)\longmapsto
 \bigl((g,B(g)),(e,B(g)^{-1}a)\bigr).
\]
Thus, the graph is an exact complement to $K$.  Conversely, the first projection from an exact complement to $G\cong E/K$ is an isomorphism of $\bfk$-schemes. The complement is the
graph of a unique morphism, and closure under multiplication gives
\eqref{eq:group-rb}.
\end{proof}

The formulation of exact complements also leads to a natural $G$-action on the space of Rota--Baxter operators.

\begin{theorem}\label{thm:gauge-all}
Let $B\in\RB_G(S)$ and $a\in G(S)$. The morphism $T_a^B$ defined by
\begin{equation}\label{eq:T}
  T_a^B(h)=a h\Ad_{B(h)}(a^{-1}),\quad \forall h\in G(S)
\end{equation}
is an $S$-isomorphism. There is a unique Rota--Baxter operator $a\cdot B$ such
that
\begin{equation}\label{eq:gauge-B}
  (a\cdot B)(T_a^B(h))=B(h).
\end{equation}
These maps define a left action of $G$ on the functor $\RB_G$, which we call the \emph{gauge action}.
\end{theorem}

\begin{proof}
We first prove a general fact about exact factorizations. If $E=HK$ is exact, then $sHs^{-1}$ is again an exact complement to $K$ for every $s\in E$. To prove this, note that one also has a exact factorization $E=KH$ by taking the inverse. Write $s=kh$.
Then $sHs^{-1}=kHk^{-1}$. The identity
\[
  (kh_0k^{-1})k_0=k h_0(k^{-1}k_0)
\]
reduces the existence and the uniqueness of the factorization $E=sHs^{-1}K$ to the original factorization. Note that the above calculation applies over $S$.

Taking the conjugation of the graph by $(a,e)$ in $G\rtimes_{\Ad}G$, we have
\begin{equation}\label{eq:graph-conj}
 (a,e)(h,B(h))(a^{-1},e)
 =\bigl(T_a^B(h),B(h)\bigr).
\end{equation}
By the general fact above, the conjugation of the graph is an exact complement of $K$, whose first projection is an isomorphism. Hence, $T_a^B$ is an isomorphism, and
\eqref{eq:gauge-B} defines the unique map, whose graph is the right-hand side of \eqref{eq:graph-conj}. It is a Rota--Baxter operator by \Cref{prop:graph}. Finally, conjugation by
$(a,e)(b,e)=(ab,e)$ proves that $a\cdot B$ defines a left action.
\end{proof}

\begin{remark}
There is another well-known elementary operation. If $B$ is a Rota--Baxter operator, then
\begin{equation}\label{eq:tilde}
\widetilde B(g):=g^{-1}B(g^{-1})
\end{equation}
is again a Rota--Baxter operator, which is known as the descent structure \cite[Proposition 2.4]{GLS}. However, the operation in \eqref{eq:tilde} is not, in
general, a gauge equivalence. We will return to this point in \Cref{ex:tilde-torus}.
\end{remark}

\subsection{Differentiation}

Let $\mathfrak g=\Lie(G)$. A linear map $R:\mathfrak g\to\mathfrak g$ is a
\emph{Rota--Baxter operator of weight one} if
\begin{equation}\label{eq:lie-rb}
 [R(x),R(y)]
 =R\bigl([R(x),y]+[x,R(y)]+[x,y]\bigr).
\end{equation}

We next describe the infinitesimal form of the action in \Cref{thm:gauge-all}. For a $\bfk$-scheme $S$ with the structure sheaf $\cO_S$, put
$\mathfrak g_S=\mathfrak g\otimes_\bfk\cO_S$, and let
$\RB_{\mathfrak g}(S)$ be the set of $\cO_S$-linear endomorphisms of
$\mathfrak g_S$ which satisfy \eqref{eq:lie-rb} on local sections.

For $a\in G(S)$ and $R\in\RB_{\mathfrak g}(S)$, set
\begin{equation}\label{eq:lie-gauge-D}
 D_a^R=\Ad(a)\circ(\id_{\mathfrak g_S}+R)-R.
\end{equation}

\begin{proposition}\label{prop:lie-gauge}
Suppose that $D_a^R$ is an automorphism of the vector bundle
$\mathfrak g_S$.  Then
\begin{equation}\label{eq:lie-gauge-action}
 a\star R:=R\circ(D_a^R)^{-1}
\end{equation}
is a weight-one Rota--Baxter operator. These formulas define a partial left
action of $G$ on $\RB_{\mathfrak g}$. More precisely, denote by $R_b=b\star R$. If $D_b^R$ is invertible, then
\begin{equation}\label{eq:lie-gauge-cocycle}
 D_a^{R_b}\circ D_b^R=D_{ab}^R.
\end{equation}

Let $\RB_{\mathfrak g}^{\mathrm{adm}}(S)$ denote the set of those $R$ for which $D_a^R$ is invertible after every base change $S'\to S$ and for every $a\in G(S')$. Then \eqref{eq:lie-gauge-action} defines a genuine action of the algebraic group $G$ on $\RB_{\mathfrak g}^{\mathrm{adm}}$ (both as functors on the category of $\bfk$-schemes). If $B\in\RB_G(S)$ and $R=d_eB$, then $R$ is admissible in the above sense, and
\begin{equation}\label{eq:diff-gauge-equivariant}
 d_eT_a^B=D_a^R,
 \quad
 d_e(a\cdot B)=a\star R.
\end{equation}
Thus, the differentiation map is equivariant with respect to the gauge actions on both sides for every algebraizable Rota--Baxter operator.
\end{proposition}

\begin{proof}
Consider the semidirect-product Lie algebra
\[
 \mathfrak e_S=\mathfrak g_S\rtimes_{\operatorname{ad}}\mathfrak g_S.
\]
The graph of $R$ is a Lie subalgebra of $\mathfrak e_S$, exactly when $R$ satisfies
\eqref{eq:lie-rb}. By differentiating the inner automorphism $\Ad_{(a,e)}$ of $G_S\rtimes_{\Ad}G_S$, we have
\[
 \Ad({(a,e)})(x,Rx)=\bigl(D_a^Rx,Rx\bigr).
\]
If $D_a^R$ is invertible, the image of the graph of $R$ is therefore the graph of $R\circ (D_a^R)^{-1}$. This proves the first assertion.

Put $R_b=R\circ (D_b^R)^{-1}$. The identities
\[
 R_b\circ D_b^R=R,\]
\[
 (\id+R_b)\circ D_b^R=\Ad(b)\circ(\id+R)
\]
imply
\[
 D_a^{R_b}\circ D_b^R
 =\Ad(a)\circ(\id+R_b)\circ D_b^R-R_b\circ D_b^R
 =\Ad(ab)\circ (\id+R)-R=D_{ab}^R.
\]
Hence, \eqref{eq:lie-gauge-action} defines a partial action on $\RB_{\mathfrak g}$ and an action on the admissible subfunctor after an arbitrary base change.

Finally, write
\[
 T_a^B(h)=a hB(h)a^{-1}B(h)^{-1}.
\]
Its relative differential at the identity is
\[
 d_eT_a^B(x)=\Ad(a)\circ(x+R(x))-R(x)=D_a^R(x).
\]
Let $S'\to S$ be a base change map and let $a'\in G(S')$. Applying
\Cref{thm:gauge-all} to $B_{S'}$ and $a'$ shows that
$T_{a'}^{B_{S'}}$ is an isomorphism. Hence, its relative differential
$D_{a'}^{R_{S'}}$ is invertible, which proves the admissibility. By differentiating
$(a\cdot B)\circ T_a^B=B$ over $S$, we get
\eqref{eq:diff-gauge-equivariant}.
\end{proof}

\begin{remark}
There is also a natural $G$-action given by changes of coordinates
\[
a\diamond R=\Ad_a\circ R\circ\Ad_{a^{-1}}.
\]
However, this action is different from the one obtained by differentiating the gauge action in \Cref{thm:gauge-all}.
\end{remark}

\begin{remark}\label{rem:lie-gauge-partial}
The admissibility condition cannot be omitted.  Let $G=\operatorname{SL}_2$ and let $\mathfrak{g}=\mathfrak{sl}_2$ with the standard basis $\{e,h,f\}$. Write
\[
 \mathfrak{sl}_2=\mathfrak{b}\oplus\mathfrak{n}^-:=(\bfk h\oplus \bfk e)\oplus \bfk f.
\]
Note that both summands are Lie subalgebras and $R=-\pr_{\bfk f}$ is a weight-one
Rota--Baxter operator, where $\pr_{\bfk f}$ denotes the projection map to $\bfk f$. For
\[
 w=\begin{pmatrix}0&-1\\1&0\end{pmatrix},
\]
one has $\Ad_w(e)=-f$.  Since $(\id+R)(e+f)=e$ and $R(e+f)=-f$, we have
\[
 D_w^R(e+f)=0.
\]
Thus, the differentiated gauge transformation does not define a global $G$-action on the whole space $\RB_{\mathfrak g}$.  In particular, this $R$ cannot be the
differential of a Rota--Baxter operator on $\operatorname{SL}_2$.
\end{remark}

\begin{proposition}\label{prop:differentiate}
If $B:G\to G$ is a Rota--Baxter operator, then
$R=d_eB$ satisfies \eqref{eq:lie-rb}.
\end{proposition}

\begin{proof}
By \Cref{prop:graph}, $\Gamma_B$ is a subgroup of $E$.  Its Lie
algebra is
the graph of $R$ in
\[
  \mathfrak g\rtimes_{\operatorname{ad}}\mathfrak g.
\]
The bracket of two elements $(x,Rx)$ and $(y,Ry)$ in the graph is given by
\[
 [(x,Rx),(y,Ry)]
 =\bigl([x,y]+[Rx,y]+[x,Ry],[Rx,Ry]\bigr),
\]
which lies in the graph of $R$ exactly when \eqref{eq:lie-rb} holds.
\end{proof}

\section{Universal enveloping algebras and algebraization}\label{sec:hopf}
We will pass from a Rota--Baxter operator on a Lie algebra to a regular operator on a specified algebraic group. The procedure consists of two stages: formal integration through the universal enveloping algebra, and an algebraization criterion in terms of its coordinate ring. We also compute the criterion explicitly for split tori.

\subsection{Rota--Baxter operators on Hopf algebras}

Let $(H,\mu,\eta,\Delta,\varepsilon,S)$ be a cocommutative Hopf algebra. We use Sweedler notation, with summation suppressed.

\begin{definition}\label{def:hopf-rb}
A \emph{Rota--Baxter operator} on $H$ is a counital coalgebra homomorphism
$\bB:H\to H$ such that
\begin{equation}\label{eq:hopf-rb}
 \bB(x)\bB(y)=
 \bB\bigl(x_{(1)}\bB(x_{(2)})yS(\bB(x_{(3)}))\bigr)
\end{equation}
for all $x,y\in H$.
\end{definition}

The antipode map provides an example of  a Rota--Baxter operator on $H$ when $H$ is cocommutative.

\begin{theorem}[{\cite[Theorem 2]{Goncharov}}]\label{thm:goncharov}
There is a bijection between the set of
Rota--Baxter operators on $U(\mathfrak g)$ and the set of Rota--Baxter operators of weight one on $\mathfrak g$, sending a Rota--Baxter operator to its restriction to the subspace of primitive elements. The extension $\bB_R$ of $R$ is characterized by
\begin{equation}\label{eq:recursion}
 \bB_R(1)=1,
 \quad
 \bB_R(xu)=R(x)\bB_R(u)-\bB_R([R(x),u])
\end{equation}
for $x\in\mathfrak g$ and $u\in U(\mathfrak g)$.
\end{theorem}

The proof first checks that the recursion is well defined on the universal enveloping algebra and then shows that the extension is a coalgebra homomorphism satisfying \eqref{eq:hopf-rb}.

Let $F_\bullet U(\mathfrak g)$ be the PBW filtration.  The recursion gives
\begin{equation}\label{eq:filtered}
  \bB_R(F_nU(\mathfrak g))\subset F_nU(\mathfrak g).
\end{equation}
Since $[\mathfrak g,F_n]\subset F_n$, we apply induction on $n$ to \eqref{eq:recursion}. This elementary observation is important for taking
the continuous dual.

\subsection{Universal enveloping algebras, distributions and duals}

Put $A=\cO(G)$ for the $\bfk$-algebra of regular functions on $G$ and let $\mathfrak m_e$ be the maximal ideal corresponding to $e$. Let $\Dist_{e,\leq n}(G)$ denote the space of
distributions of order at most $n$ supported at $e$. Moreover, we have
\[
 \Dist_{e,\leq n}(G)=(A/\mathfrak m_e^{n+1})^*.
\]
Their union $\Dist_e(G)$ is a filtered Hopf algebra with the multiplication being the convolution against the pullback by the group multiplication \cite[I.7.7]{Jantzen}. Our convention identifies $\mathfrak g$ with the subspace of primitive distributions and gives its usual Lie bracket.

By Cartier's theorem \cite[I.7.10]{Jantzen}, in characteristic zero, the canonical map
\begin{equation}\label{eq:cartier}
 U(\mathfrak g)\xrightarrow{\sim}\Dist_e(G)
\end{equation}
is an isomorphism of filtered Hopf algebras. Moreover, there is an identification between the associated graded algebra of $U(\mathfrak g)$
and $\operatorname{Sym}(\mathfrak g)$ by the PBW theorem.

The perfect pairings on finite jets imply there is a topological isomorphism
\begin{equation}\label{eq:completion-dual}
 \widehat{\cO}_{G,e}
 \simeq \varprojlim_n (F_nU(\mathfrak g))^*.
\end{equation}
We write
\[
  \widehat G_e=\operatorname{Spf}(\widehat{\cO}_{G,e})
\]
for the formal completion of $G$ at $e$. For a linear map
$T:U(\mathfrak g)\to U(\mathfrak g)$, its transpose is defined by
\[
  T^\vee(\varphi)=\varphi\circ T.
\]
The inverse limit on the right-hand side of \eqref{eq:completion-dual} is naturally the full linear dual $U(\mathfrak g)^*$, endowed with its inverse-limit topology. Here, the finite dimensional duals of the jet spaces $(F_nU(\mathfrak g))^*$ agree with $\cO(G)/\mathfrak m_e^{n+1}$ through the canonical pairing between regular functions and distributions.

There is another notion of dual for $U(\mathfrak g)$ known as \emph{the finite dual} $U(\mathfrak g)^\circ$, which consists of linear functionals with kernels containing a two-sided ideal
of finite codimension. Equivalently, $U(\mathfrak g)^\circ$ is spanned by matrix coefficients of finite dimensional $U(\mathfrak g)$-modules.  

Note that there is a distinguished map
\begin{equation}\label{eq:iota}
\begin{aligned}
\iota_G:A&\rightarrow U(\mathfrak g)^*,\\
f&\mapsto (u\mapsto u(f)).
 \end{aligned}
\end{equation}

\begin{lemma}\label{lem:iota}
Assume that $G$ is connected. Then the map $\iota_G$ is an injective algebra homomorphism. Its image is contained in the finite dual $U(\mathfrak g)^\circ$.
\end{lemma}

\begin{proof}
Although it seems known to experts, we still provide a short proof for completeness. The coalgebra structure on $U(\mathfrak g)$ makes \eqref{eq:iota} an algebra homomorphism. Moreover, every regular function belongs to a finite dimensional subcomodule of
the regular representation of $G$. After the differentiation,
$u\mapsto u(f)$ is a matrix coefficient of a finite dimensional
$U(\mathfrak g)$-module. Hence, it belongs to the finite dual.

If $\iota_G(f)=0$, then every distribution supported at $e$ annihilates $f$.
By the perfectness of the jet pairing and the Krull intersection theorem, we know that the germ of $f$ at $e$ is zero. Since a connected smooth algebraic group is irreducible, its coordinate ring $A$ is a domain. Hence, the map $A\to\cO_{G,e}$ is injective. It follows that $f=0$.
\end{proof}

\subsection{Formal integration and a criterion for algebraization}

Since $\bB_R$ is a filtered coalgebra homomorphism, its transpose is a continuous local algebra endomorphism of $\widehat{\cO}_{G,e}$. Let
\[
  \widehat\beta_R:\widehatG\longrightarrow\widehatG
\]
be the corresponding pointed formal map.

Here, a \emph{formal Rota--Baxter operator} means a pointed formal self-map of $\widehat G_e$, for which \eqref{eq:group-rb} holds on the formal completion of $G\times G$ at $(e,e)$.

The passage from the transpose to a map on $G$ involves the duality between affine schemes and their coordinate rings. When the distinguished copy of $A=\cO(G)$ is stable under the transported map, the map is a unital algebra endomorphism on $A$, which is then the pullback of a regular map $G\to G$. It might not be a Hopf algebra endomorphism, since a Rota--Baxter operator need not be a group homomorphism.  

\begin{theorem}\label{thm:algebraization}
Assume that $G$ is connected. The formal map $\widehat\beta_R$ is a formal Rota--Baxter operator and has the tangent map $R$. It is the completion of a unique pointed regular map
$\beta_R:G\to G$ if and only if
\begin{equation}\label{eq:stability}
 \bB_R^\vee\bigl(\iota_G(A)\bigr)\subseteq\iota_G(A).
\end{equation}
When \eqref{eq:stability} holds, $\beta_R$ is a Rota--Baxter operator on $G$ with $d_e\beta_R=R$.
\end{theorem}

\begin{proof}
Note that the transpose of a counital coalgebra homomorphism is a unital algebra homomorphism. Since $\bB_R^\vee$ preserves the filtration in \eqref{eq:filtered}, we know that it is continuous by \eqref{eq:completion-dual}. It is local because $\bB_R(1)=1$. Hence, $\widehat\beta_R$ is a pointed formal map. On the tangent space, its pullback is
dual to the restriction of $\bB_R$ to the subspace of primitive distributions. Hence, the tangent map of $\widehat\beta_R$ is $R$.

We first verify the formal Rota--Baxter identity for $\widehat\beta_R$. Take
\[
 F_L(g,h):=\widehat\beta_R(g)\widehat\beta_R(h),
\]
\[
 F_R(g,h):=\widehat\beta_R
 \bigl(g\widehat\beta_R(g)h\widehat\beta_R(g)^{-1}\bigr).
\]
The distributions supported at $(e,e)$ are spanned by pure tensors $x\otimes y$. By the functoriality of distributions, we have
\[
 (F_L)_*(x\otimes y)=\bB_R(x)\bB_R(y),
\]
\[
 (F_R)_*(x\otimes y)=
 \bB_R\bigl(x_{(1)}\bB_R(x_{(2)})y
 S(\bB_R(x_{(3)}))\bigr).
\]
These two expressions agree by \eqref{eq:hopf-rb}. The pairings on finite jets are perfect so that the formal maps agree.

Suppose \eqref{eq:stability} holds. By \Cref{lem:iota}, there is a unique algebra homomorphism $\beta_R^*:A\to A$ such that
\[
 \iota_G\circ\beta_R^*=\bB_R^\vee\circ\iota_G.
\]
By the duality between affine schemes and their coordinate rings, $\beta_R^*$ is the pullback of a unique regular map $\beta_R:G\to G$. Note that the map is pointed, since evaluation at $e$ is exactly given by pairing with $1\in U(\mathfrak g)$. For $u\in\Dist_e(G)$ and $f\in A$,
we have
\[
 \langle(\beta_R)_*u,f\rangle
 =\langle u,\beta_R^*f\rangle
 =\langle\bB_R(u),f\rangle.
\]
Hence, $(\beta_R)_*=\bB_R$ on $\Dist_e(G)$.  The converse follows from the
same equality.

It remains to pass from the formal identity to a global one. For $f\in A$, the difference $F_L^*f-F_R^*f$ is zero in the completed local ring $\widehat{\cO}_{G\times G,(e,e)}$. Note that the canonical map from a Noetherian local ring to its completion is injective by the Krull intersection theorem. Since $G\times G$ is irreducible, the map from $\cO(G\times G)$ to its localization $\cO(G\times G)_{(e,e)}$ is also injective. Hence, a regular function vanishes in $\widehat{\cO}_{G\times G,(e,e)}$ implies that it vanishes globally. Thus, the two regular maps agree on $G\times G$, i.e., we have
\[
 \beta_R(g)\beta_R(h)
 =\beta_R\bigl(g\Ad_{\beta_R(g)}(h)\bigr),
\]
which is precisely the identity in \Cref{def:group-rb}. By taking
the restriction of $(\beta_R)_*=\bB_R$ to the space of primitive distributions of degree one, we obtain $d_e\beta_R=R$.
\end{proof}

\begin{corollary}\label{cor:differentiation-injective}
Let $G$ be connected. The differentiation map defines an injective map
\[
d_e:\RB_G(\bfk)\longrightarrow\RB_{\mathfrak g}(\bfk).
\]
Its image consists exactly of the operators $R$ satisfying
\eqref{eq:stability}. For such an $R$, the unique Rota--Baxter operator whose differential is $R$ is $\beta_R$ in \Cref{thm:algebraization}.
\end{corollary}

\begin{proof}
Let $B\in\RB_G(\bfk)$ and put $R=d_eB$. The pushforward map on the space of distributions induces a counital coalgebra homomorphism
\[
B_*:\Dist_e(G)\longrightarrow\Dist_e(G).
\]
Adapting the proof of
\Cref{thm:algebraization} for $(F_L)_*$ and $(F_R)_*$, we find that $B_*$ satisfies \eqref{eq:hopf-rb}, and its restriction to the subspace of primitive distributions is $R$. Then the uniqueness in \Cref{thm:goncharov} implies
$B_*=\bB_R$. Consequently, we have
\[
 \iota_G\circ B^*=\bB_R^\vee\circ\iota_G.
\]
Thus, \eqref{eq:stability} holds and $B=\beta_R$. The converse is just \Cref{thm:algebraization}.
\end{proof}

\begin{proposition}\label{prop:isomorphism-functoriality}
Let $\varphi:G\xrightarrow{\sim}G'$ be an isomorphism of connected algebraic groups, and let $R\in\RB_{\mathfrak g}(\bfk)$.  Put
\[
 R':=d_e\varphi\circ R\circ(d_e\varphi)^{-1}.
\]
Then $R$ algebraizes on $G$ if and only if $R'$ algebraizes on $G'$.  In
this case,
\[
 \beta_{R'}=\varphi\circ\beta_R\circ\varphi^{-1}.
\]
\end{proposition}

\begin{proof}
Taking conjugation by $d_e\varphi$ converts the Rota--Baxter identity for $R$ to the identity for $R'$. Suppose that $R$ algebraizes. Then the map $\varphi\beta_R\varphi^{-1}$ is a Rota--Baxter operator on $G'$ and its differential is $R'$.  By \Cref{cor:differentiation-injective}, it is identified with $\beta_{R'}$. Applying the same argument to
$\varphi^{-1}$ proves the converse.
\end{proof}

\begin{remark}\label{rem:differential-versus-jets}
Note that an arbitrary pointed regular map is not fully determined by its differential even under the assumption that the algebraic group is connected. For example, $x\mapsto x$ and $x\mapsto x+x^2$ on $\Ga$ have the same differential at the origin. In \Cref{cor:differentiation-injective}, the Rota--Baxter identity and \Cref{thm:goncharov} identify the maps on $\Dist_e(G)$, not just on $\mathfrak g$. Then, the injectivity of $\iota_G$ determines the global regular maps.
\end{remark}

\subsection{The criterion for split tori}

\begin{theorem}\label{thm:torus-algebraization}
Let $T$ be a split torus and let $X_*(T)$ be the cocharacter lattice. Denote $\mathfrak t=\Lie(T)$. Under the canonical identification
\[
 \mathfrak t=X_*(T)\otimes_{\mathbb Z}\bfk,
\]
every $R\in\End_\bfk(\mathfrak t)$ is a Rota--Baxter operator on $\mathfrak t$.  The formal map $\widehat\beta_R$ algebraizes on $T$ if and only if
\begin{equation}\label{eq:cocharacter-integrality}
 R\bigl(X_*(T)\bigr)\subseteq X_*(T).
\end{equation}
When this holds, $\beta_R$ is the endomorphism of $T$ induced by the integral endomorphism $R|_{X_*(T)}$. 
\end{theorem}
In terms of a basis of $X_*(T)$, the condition in \eqref{eq:cocharacter-integrality} says that the matrix corresponding to $R$ has integral entries.

\begin{proof}
Since $\mathfrak t$ is abelian, the Rota--Baxter identity holds for every linear map $R$. Suppose first that $R$ algebraizes. Since $T$ is commutative, by \Cref{def:group-rb}, $\beta_R$ is an endomorphism of $T$. Its action on the cocharacter lattice is an integral endomorphism of $X_*(T)$, whose scalar extension to $\mathfrak t$ coincides with the differentiation of $\beta_R$. By \Cref{thm:algebraization}, we have $d_e\beta_R=R$. Hence \eqref{eq:cocharacter-integrality} holds.

Conversely, suppose that $R$ preserves $X_*(T)$. Then the restriction of $R$ to this lattice defines a unique torus endomorphism $f_R:T\to T$, which is a Rota--Baxter operator since $T$ is commutative and its differential is $R$. By \Cref{cor:differentiation-injective}, $R$ satisfies the algebraization criterion and $f_R=\beta_R$.
\end{proof}

\begin{proposition}\label{prop:torus-obstruction}
Let $G=\Gm$. Let $X$ be the map
\[
 f\longmapsto\left.t\frac{df}{dt}\right|_{t=1},\quad\forall f\in \mathcal{O}(\Gm),
\]
so that $\mathfrak g=\bfk X$. Take
$R(X)=\lambda X$. Then we have
\[
  \bB_R(X^n)=\lambda^nX^n.
\]
The formal map in \Cref{thm:algebraization} is $t\mapsto t^\lambda$. It algebraizes on $\Gm$ if and only if $\lambda\in\mathbb Z\subset \bfk$. In particular, the Rota--Baxter operator $R=\frac12\id$ does not algebraize on $\Gm$.
\end{proposition}

Here $t^\lambda$ denotes the binomial series
\[
  \bigl(1+(t-1)\bigr)^\lambda
  =\sum_{j\geq0}\binom{\lambda}{j}(t-1)^j
  \quad\text{in }\widehat{\cO}_{\Gm,1}.
\]

\begin{proof}
Since $\mathfrak g$ is abelian, the commutator term in \eqref{eq:recursion} vanishes. We get the corresponding formula for $\bB_R$ by induction. Under the natural pairing between regular functions and distributions, we have
\begin{equation}\label{eq:torus-pairing}
 \langle X^n,t^m\rangle=m^n,
 \quad \forall m\in\mathbb Z.
\end{equation}
Thus, $t^m$ is the evaluation functional $\operatorname{ev}_m$ on $\bfk[X]$ and $\bB_R^\vee(\operatorname{ev}_a)=\operatorname{ev}_{\lambda a}$. In the formal coordinate $s=\log t$, the transpose sends $f(s)$ to
$f(\lambda s)$. This is the pullback of the formal map $t\mapsto t^\lambda$. Since the cocharacter lattice of $\Gm$ is $\mathbb Z X$, the last two assertions follow from \Cref{thm:torus-algebraization}.
\end{proof}

The element $\operatorname{ev}_{1/2}$ lies in the finite dual because its kernel contains the ideal $(X-1/2)$. However, it does not lie in the copy of $\bfk[t,t^{-1}]$. This indeed shows that
\[
  \cO(\Gm)\subsetneq U(\bfk X)^\circ.
\]
Furthermore, the full dual is even larger and describes the completed formal object.

\section{1-Cocycles and the coboundary action}\label{sec:cocycles}

We compare Rota--Baxter operators with nonabelian regular conjugation 1-cocycles and group endomorphisms. This dictionary identifies the invertibility condition and transports gauge
equivalence into an algebraic conjugation action, which also works for nonreduced parameter schemes. Indeed, this dictionary also motivates our discovery of the action in \Cref{thm:gauge-all}.

\subsection{The dictionary}\label{sec:dic}

A \emph{regular conjugation 1-cocycle} is a morphism $c:G\to G$ such that
\begin{equation}\label{eq:cocycle}
  c(xy)=c(x)\Ad_x(c(y)),\quad \forall x,y\in G.
\end{equation}
The left coboundary $G$-action on the set of regular conjugation 1-cocycles is defined by
\begin{equation}\label{eq:coboundary}
  (a\star c)(x):=a c(x)\Ad_x(a^{-1}),\quad\forall a,x\in G.
\end{equation}
It is direct to check that $a\star c$ is also a regular conjugation 1-cocycle. It satisfies $a\star(b\star c)=(ab)\star c$ for any $a,b\in G$.

We use the notation
\[
 H^1_{\mathrm{alg}}(G,G_{\Ad})
:=Z^1_{\mathrm{alg}}(G,G_{\Ad})(\bfk)/G(\bfk)
\]
for this orbit set. We remark that this notation does not mean fppf cohomology of $\Spec \bfk$. It also does not include arbitrary set-theoretic 1-cocycles.

\begin{proposition}\label{prop:cocycle-end}
The assignments
\begin{equation}\label{eq:cocycle-end}
 \theta_c(x)=c(x)x,
 \quad
 c_\theta(x)=\theta(x)x^{-1}
\end{equation}
define an isomorphism and its inverse between the functor of regular conjugation
1-cocycles and the functor $\underline{\End}_{\mathrm{gp}}(G)$. Under this isomorphism,
the action in \eqref{eq:coboundary} becomes
\[
  \theta\longmapsto\Ad_a\circ\theta.
\]
\end{proposition}

\begin{proof}
By the cocycle identity, we have
\[
 \theta_c(xy)=c(x)x c(y)y=\theta_c(x)\theta_c(y),\quad \forall x,y\in G.
\]
The converse follows from similar calculation in the reverse direction. Note that we also have
\[
 \theta_{a\star c}(x)
 =a c(x)x a^{-1}=\Ad_a(\theta_c(x)).
\]
All formulas commute with base change.
\end{proof}

Let $Z^1_{\mathrm{alg,iso}}(G,G_{\Ad})$ be the subfunctor consisting of 1-cocycles which
are isomorphisms of schemes.

\begin{theorem}\label{thm:inverse-cocycle}
There is an isomorphism of functors
\begin{equation*}
\begin{aligned}
\RB_G^\times&\xrightarrow{\sim}
 Z^1_{\mathrm{alg,iso}}(G,G_{\Ad}),\\
 B&\longmapsto B^{-1},
\end{aligned}
\end{equation*}
where $\RB_G^\times$ denotes the functor parameterizing Rota--Baxter operators that are isomorphisms of schemes and $B^{-1}$ denotes the functional inverse. The isomorphism intertwines the gauge action given in \Cref{thm:gauge-all} with the coboundary action
\eqref{eq:coboundary}.
\end{theorem}

\begin{proof}
Put $c=B^{-1}$. In \eqref{eq:group-rb}, substitute $g=c(x)$ and $h=c(y)$.
Applying $c$, we have
\[
 c(xy)=c(x)\Ad_x(c(y)).
\]
We can reverse the above calculation to show that every isomorphic 1-cocycle corresponds to an invertible Rota--Baxter operator.

To prove the correspondence between these two actions, set $h=c(x)$ in \eqref{eq:gauge-B}. Since $B(h)=x$, by combining \eqref{eq:coboundary},we
obtain
\[
 (a\cdot B)^{-1}(x)=a c(x)\Ad_x(a^{-1})=(a\star B^{-1})(x).
\]
\end{proof}

Define a subfunctor $\DEnd(G)$ of $\underline{\End}_{\mathrm{gp}}(G)$ by
\begin{equation}\label{eq:Dfunctor}
 \DEnd(G)(S)=
 \left\{\theta\in\End_{S\text{-}\mathrm{gp}}(G_S):
 \delta_\theta(x)=\theta(x)x^{-1}
 \text{ is an $S$-isomorphism}\right\}.
\end{equation}
By \Cref{prop:cocycle-end} and \Cref{thm:inverse-cocycle}, this is the functor parameterizing invertible Rota--Baxter operators.

Indeed, the same difference construction works for Lie algebras, which does not require an algebraization hypothesis.

\begin{proposition}\label{prop:lie-difference}
Let $\mathfrak g$ be a finite-dimensional Lie algebra. The assignments
\begin{equation}\label{eq:lie-difference}
\begin{aligned}
 R&\longmapsto\alpha_R=\id+R^{-1},\\
 \alpha&\longmapsto R_\alpha=(\alpha-\id)^{-1}
\end{aligned}
\end{equation}
define inverse bijections between
\[
 \RB_{\mathfrak g}(\bfk)\cap\operatorname{GL}(\mathfrak g)
\]
and the set of Lie algebra endomorphisms $\alpha$, for which
$\alpha-\id$ is invertible.

If $B$ is an invertible Rota--Baxter operator on $G$, let
\[
 \theta_B(x)=B^{-1}(x)x,
\]
where $B^{-1}$ denotes the functional inverse.  Then
\begin{equation}\label{eq:bridge-differential}
 d_e\theta_B=\id+(d_eB)^{-1}.
\end{equation}
Hence, we have the following commutative diagram. Moreover, both horizontal arrows are bijections:
\begin{equation}\label{eq:bridge-diagram}
\begin{array}{ccc}
 \RB_G^\times(\bfk) & \xrightarrow{\ B\mapsto\theta_B\ } & \DEnd(G)(\bfk)\\[4pt]
 \big\downarrow d_e && \big\downarrow d_e\\[4pt]
 \RB_{\mathfrak g}(\bfk)\cap\operatorname{GL}(\mathfrak g)
 & \xrightarrow{\ R\mapsto\id+R^{-1}\ } &
 \{\alpha\in\End_{\mathrm{Lie}}(\mathfrak g):
       \alpha-\id\in\operatorname{GL}(\mathfrak g)\}.
\end{array}
\end{equation}
\end{proposition}

\begin{proof}
Suppose first that $R$ is invertible and put $D=R^{-1}$. In the Rota--Baxter identity on $\mathfrak g$, substitute $x=D(a)$ and $y=D(b)$, and then apply $D$. We get
\begin{equation}\label{eq:D-identity}
 D[a,b]=[D(a),b]+[a,D(b)]+[D(a),D(b)].
\end{equation}
Consequently, we have 
\[
 [a+D(a),b+D(b)]=[a,b]+D[a,b],
\]
which means that $\id+D$ is a Lie algebra endomorphism. Note that the above calculation can be reversed. Hence, an endomorphism $\alpha$ with $\alpha-\id$ invertible gives the Rota--Baxter operator $(\alpha-\id)^{-1}$.

For the statement for groups, note that \Cref{thm:inverse-cocycle} and \Cref{prop:cocycle-end} give the upper horizontal bijection. The differential of $B^{-1}$ is $(d_eB)^{-1}$. By differentiating
$\theta_B(x)=B^{-1}(x)x$, we have \eqref{eq:bridge-differential}. This also proves the commutativity of \eqref{eq:bridge-diagram}.
\end{proof}

\subsection{Comparison with algebraic nonabelian cohomology}

\begin{example}\label{ex:h1-torus}
Take $G=\Gm$. Then a 1-cocycle is simply a group homomorphism. Its pullback sends $t$ to a unit of $\bfk[t,t^{-1}]$. Every such unit has the form $a t^n$, and the identity condition gives $a=1$. Thus, the regular 1-cocycles are exactly
\[
  c_n(t)=t^n,\quad\forall n\in\mathbb Z.
\]
The map $c_n$ is an automorphism exactly when $n=\pm1$.  Hence, we have
\[
 H^1_{\mathrm{alg}}(\Gm,\Gm)=\mathbb Z,
 \quad
 Z^1_{\mathrm{alg,iso}}(\Gm,\Gm)/\Gm=\{1,-1\}.
\]
Thus, the locus consisting of isomorphic 1-cocycles is a proper subset of algebraic
nonabelian $H^1$. We remark that the isomorphism condition is a scheme-theoretic condition instead of a set-theoretic condition. For example, in characteristic $p$ (although beyond the scope of this article), the map $t\mapsto t^p$ is bijective on geometric points but is not an isomorphism of schemes.
\end{example}

We write
\begin{equation}\label{eq:h1diff}
 H^1_{\mathrm{diff}}(G,G_{\Ad})
 =\DEnd(G)(\bfk)/G(\bfk).
\end{equation}
The natural map from this set to $H^1_{\mathrm{alg}}(G,G_{\Ad})$ is
injective, and its image consists exactly of the classes having an
isomorphic 1-cocycle representative. The notation in \eqref{eq:h1diff} reflects the condition in \eqref{eq:Dfunctor}.
\begin{remark}\label{rem:two-conditions}
The condition in \Cref{thm:algebraization} concerns the existence of a regular map with a prescribed formal completion. The condition defining $H^1_{\mathrm{diff}}$ concerns the invertibility of a 1-cocycle attached to an algebraic map. \Cref{prop:torus-obstruction} constructs a formal map that fails the first condition, while \Cref{ex:h1-torus} gives algebraic 1-cocycles that fail
the second. 
\end{remark}

\section{Moduli of invertible Rota--Baxter operators for reductive groups}\label{sec:moduli}

For connected reductive groups, rigidity of the
endomorphism scheme turns the classification problem of invertible
Rota--Baxter operators into a particularly transparent moduli
problem. We first recall the rigidity result in \cite{Brion} and then describe the parameter scheme, quotient stack, and the corresponding good moduli space.

\subsection{Rigidity of the endomorphism scheme}

An affine algebraic group $H$ is \emph{linearly reductive} if the functor of taking invariants $V\mapsto V^H$ on the category of rational $H$-modules is exact. We now assume that
$G$ is connected and reductive. In characteristic zero, this implies that $G$ is linearly reductive \cite[Remark 12.3]{Alper}.

Let $E_G=\underline{\Hom}_{\mathrm{gp}}(G,G)$ denote the endomorphism scheme with 
\[
E_G(S)=\Hom_{S\text{-}\mathrm{gp}}(G_S,G_S).
\]
The group $G$ acts on $E_G$ by conjugation on the target, i.e.,
\[a\cdot\theta:=\Ad_a\circ\theta,\quad\forall a\in G,\theta\in E_G.\]
The following theorem provides the
rigidity result which will be used later.

\begin{theorem}[{\cite{Brion}}]\label{thm:brion}
The functor $E_G$ is represented by a smooth scheme locally of finite type.
Every connected component is affine of finite type. Each connected component is a single scheme-theoretic orbit for the conjugation action on the target.
\end{theorem}

Representability, smoothness, and openness of the orbits follow from
\cite[Theorem 2]{Brion}. Moreover, by \cite[Corollary 4.9]{Brion}, every connected component is affine of finite type. Furthermore, by \cite[Remark 6.4]{Brion}, connected components of the representing scheme are exactly the orbits of the identity component under the conjugation action on the target. Since $G$ is connected, these components are exactly the above orbits for the conjugation action on the target.

\subsection{The parameter scheme and quotient stack}

Let $D_G$ denote the subfunctor \eqref{eq:Dfunctor}.  For a $G$-scheme $X$,
the quotient stack $[X/G]$ assigns to a $\bfk$-scheme $S$ the groupoid of pairs
$(P,f)$, where $P\to S$ is a $G$-torsor for the fppf topology and
$f:P\to X$ is a $G$-equivariant morphism. For an algebraic group $C$, write
$BC=[\Spec \bfk/C]$ for the classifying stack.

A morphism $q:\mathcal X\to M$ from an algebraic stack to an algebraic space is a \emph{good moduli space} if it is quasi-compact, the natural map $\cO_M\to q_*\cO_{\mathcal X}$ is an isomorphism, and $q_*$ is exact on the category of quasi-coherent sheaves
\cite[Definition 4.1]{Alper}. We have the following result for moduli.

\begin{theorem}\label{thm:main-moduli}
Let $G$ be a connected reductive algebraic group over $\bfk$.
\begin{enumerate}[label=\textup{(\roman*)}]
\item The functor $D_G$ is represented by a smooth open and closed subscheme
of $E_G$. It is a disjoint union of affine components of finite type.
\item Let $I_G$ be the set of these components.  Choose
$\theta_i\in D_G(\bfk)$ in each component and put
\[
  C_i=Z_G(\theta_i(G)).
\]
Then every $C_i$ is linearly reductive and there is a decomposition
\begin{equation}\label{eq:stack-decomposition}
 \cM_{\RB}^\times(G):=[D_G/G]
 \simeq\coprod_{i\in I_G}BC_i.
\end{equation}
\item The morphism
\begin{equation}\label{eq:good-moduli}
 q:\cM_{\RB}^\times(G)
 \longrightarrow M_{\RB}^\times(G):=
 \coprod_{i\in I_G}\Spec \bfk
\end{equation}
is a good moduli space.
\item For every algebraically closed field extension $\Omega/\bfk$, the points of
$M_{\RB}^\times(G)_\Omega$ are in natural bijection with the isomorphism
classes in $\cM_{\RB}^\times(G)(\Omega)$.  They are also in bijection with
\[
 H^1_{\mathrm{diff}}(G_\Omega,(G_\Omega)_{\Ad}).
\]
\end{enumerate}
\end{theorem}

\begin{proof}
We first show that the condition in \eqref{eq:Dfunctor} is constant on each orbit. Let $c=\delta_\theta$ and let $a\in G(S)$. The difference map associated to $\Ad_a\circ\theta$ is
\[
 (a\star c)(x)=a c(x)\Ad_x(a^{-1}).
\]
Assume that $c$ is an isomorphism. Choose $h=c^{-1}(a^{-1})$. By the cocycle identity, we have
\begin{equation}\label{eq:translation-cocycle}
 (a\star c)(x)=a c(xh).
\end{equation}
Thus, $a\star c=L_a\circ c\circ R_h$, which implies that $a\star c$ is also an isomorphism. 

Let $D$ be the union of those connected components of $E_G$ containing a point with an invertible difference map. It is both open and closed. We must check that it represents $D_G$ for arbitrary families. Fix a component and a
point $\theta_i$ in it. The orbit map
\[
 G\longrightarrow G\cdot\theta_i
\]
is smooth and surjective. Since it is also of finite type, it is an fpqc cover. Pulling back the universal difference map along this cover yields the family in \eqref{eq:translation-cocycle}, whose members are isomorphisms when $\delta_{\theta_i}$ is. The property of being an isomorphism of schemes descends fpqc-locally on the base. Thus, the universal difference map is an isomorphism on each selected component.

Conversely, assume that an $S$-point of $E_G$ has an invertible difference map. Every
geometric fiber lies in a selected component by an explicit calculation of the orbit. Its
topological image therefore lies in the open subscheme $D$. The map
$S\to E_G$ factors through $D$, including when $S$ is nonreduced. Hence, we prove (i).

Note that the stabilizer of $\theta_i$ is $C_i$.  By \Cref{thm:brion}, its orbit is
the affine scheme $G/C_i$.  The orbit map is smooth, so $C_i$ is smooth.
Since $G$ is linearly reductive and $G/C_i$ is affine, Matsushima's theorem (in Alper's form) \cite[Theorem 12.15(i)]{Alper} shows that $C_i$ is linearly reductive. The standard quotient-stack equivalence now gives
\[
 [G/C_i\,/G]\simeq BC_i.
\]
Hence, we prove (ii).

It is well known that quasi-coherent sheaves on $BC_i$ correspond to rational $C_i$-modules, and
the pushforward map along $BC_i\to\Spec \bfk$ gives the functor of taking invariants. Since $C_i$ is linear reductive, this functor is exact. Moreover, the space of invariants in the structure sheaf is $\bfk$. Hence $BC_i\to\Spec \bfk$ is a good moduli space \cite[Definition 12.1 and Remark 12.2]{Alper}. Since quasi-coherent sheaves on a disjoint union decompose componentwise, the same two properties hold for \eqref{eq:good-moduli}. This morphism is quasi-compact even when $I_G$ is infinite: a quasi-compact open subset of the target meets only finitely many components. Hence, we prove (iii).

Another crucial fact is that any torsor under an affine algebraic group over an algebraically closed field is trivial. To see this, note that such a torsor is a nonempty scheme of finite type. Hence, it has a rational point, which can be used to trivialize the torsor. Hence, isomorphism classes in $[D_G/G](\Omega)$ are exactly parameterized by $D_G(\Omega)/G(\Omega)$. By \Cref{prop:cocycle-end} and \Cref{thm:inverse-cocycle}, this is precisely the set in \eqref{eq:h1diff}. Since each component is a single orbit, each component contributes one point to the good moduli space. We finish the proof.
\end{proof}

\begin{remark}\label{rem:stack-vs-space}
We remind readers about the differences between a quotient stack and a good moduli scheme. The quotient stack in \eqref{eq:stack-decomposition} records objects together with their automorphism groups. The corresponding good moduli scheme records geometric isomorphism classes, but forgets the automorphism groups $C_i$. The description in terms of
isomorphism classes holds here because every component is a single closed
orbit. In general, a good quotient identifies the orbits whose closure meet.
\end{remark}

\section{Examples}\label{sec:examples}

We provide various examples for algebraization and invertibility. We also illustrate the different behavior for 
reductive and nonreductive groups. The tori exhibit certain integral constraints; while the simple and additive groups display contrasting rigidity and deformation phenomena.

\subsection{Split tori}

Let $T=\Gm^r$. There is a natural identification between $\Hom_{\mathrm{gp}}(T,T)$ and
$\End_{\mathbb Z}(X^*(T))=M_r(\mathbb Z)$ by sending an endomorphism of $T$ to its pullback on $X^*(T)$. Since the conjugation action on $T$ is trivial, \Cref{thm:brion} makes each homomorphism a connected component. Therefore, we have
\[
  E_T=\coprod_{A\in M_r(\mathbb Z)}\Spec \bfk.
\]
The difference map of $A$ has the corresponding matrix $A-I_r$. Therefore, we have
\begin{equation}\label{eq:torus-D}
  D_T=\coprod_{A-I_r\in\operatorname{GL}_r(\mathbb Z)}\Spec \bfk.
\end{equation}
In terms of the cocycle parameterization, the indexing set is
$\operatorname{GL}_r(\mathbb Z)$. Thus, we have
\begin{equation}\label{eq:torus-stack}
 \cM_{\RB}^\times(T)
 =\coprod_{U\in\operatorname{GL}_r(\mathbb Z)}BT,
 \quad
 M_{\RB}^\times(T)
 =\coprod_{U\in\operatorname{GL}_r(\mathbb Z)}\Spec \bfk.
\end{equation}

For $r\geq2$, the indexing set is infinite. The moduli scheme is not quasi-compact, not of finite type, and not affine. Thus, the discrete good moduli space in \Cref{thm:main-moduli} may not arise as a single affine GIT quotient of finite type.

For $T=\Gm$, the two invertible Rota--Baxter operators are the identity map and the inversion map, both of which have stabilizer $T$. Hence, the moduli stack is
\[
  BT\amalg BT,
\]
and the good moduli scheme consists of two points.

\begin{example}\label{ex:tilde-torus}
For $B_n(t)=t^n$, the descent structure in \eqref{eq:tilde} is
\[
  \widetilde B_n(t)=t^{-n-1}.
\]
The gauge action is trivial because $T$ is commutative. Hence $B_n$ and $\widetilde B_n$ are generally not gauge equivalent. In particular, the inversion map is sent to the constant map. The operation $B\mapsto\widetilde B$
does not preserve the invertible locus.
\end{example}

\subsection{Simple algebraic groups and product constructions}

\begin{proposition}\label{prop:simple}
Let $G$ be a connected simple algebraic group over $\bfk$. By being simple, we mean simple in the algebraic setting, which can be characterized by the condition that the corresponding Lie algebra $\mathfrak g$ is simple (sometimes called almost simple). The inversion map is the only invertible Rota--Baxter operator on $G$. Consequently, we have
\[
  \cM_{\RB}^\times(G)\simeq BG,
  \quad
  M_{\RB}^\times(G)\simeq\Spec \bfk.
\]
\end{proposition}

\begin{proof}
Let $B$ be an invertible Rota--Baxter operator on $G$ and put
$R=d_eB$. Then $R$ is invertible. By
\Cref{prop:lie-difference},
$\alpha=\id+R^{-1}$ is a Lie algebra endomorphism of $\mathfrak g$
with $\alpha-\id$ invertible. The fact that $\mathfrak g$ is simple implies that either $\alpha=0$ or $\alpha$ is an automorphism. In the latter case, Jacobson's fixed-point theorem \cite[Theorem 9]{Jacobson} states that there exists a non-zero element
$x\in\mathfrak g$ with $\alpha(x)=x$, contradicting
invertibility of $\alpha-\id$.  Thus $\alpha=0$, so $R=-\id$.

Note that the inversion map $g\mapsto g^{-1}$ is a Rota--Baxter operator on $G$ with the differential $-\id$. By the injectivity of the differentiation map
in \Cref{cor:differentiation-injective}, we know that $B(g)=g^{-1}$. The associated difference endomorphism is the trivial homomorphism, whose stabilizer under conjugation is the whole group $G$. Furthermore, a connected simple affine algebraic group is reductive. Thus, the statements about moduli follow from \Cref{thm:main-moduli}.
\end{proof}

Explicit constructions of Rota--Baxter operators on products of groups have appeared previously (see \cite[Proposition 4.4]{BG}). The following endomorphisms constructed by shifts provide an explicit family of invertible Rota--Baxter operators. Furthermore, we can compute the corresponding stabilizers and components explicitly.

\begin{proposition}\label{prop:shift-family}
Let $H$ be an algebraic group, let $m\geq2$, and put $G=H^m$.  Define
\[
 \theta(x_1,\ldots,x_m)=(e,x_1,\ldots,x_{m-1}).
\]
Then $\theta$ is in $\DEnd(G)(\bfk)$.  The associated invertible Rota--Baxter operator is
\begin{equation}\label{eq:shift-RB}
 B(y_1,\ldots,y_m)
 =\bigl(y_1^{-1},(y_1y_2)^{-1},\ldots,(y_1\cdots y_m)^{-1}\bigr).
\end{equation}
If $H$ is connected and reductive, then the component of the quotient stack containing $B$ is
\[
 B\bigl(H\times Z(H)^{m-1}\bigr).
\]
In particular, it is $BH$ when $H$ is adjoint.
\end{proposition}

\begin{proof}
Note that the map $\theta$ is an endomorphism of $G$, and its difference map is
\[
 \delta_\theta(x_1,\ldots,x_m)
 =\bigl(x_1^{-1},x_1x_2^{-1},\ldots,x_{m-1}x_m^{-1}\bigr).
\]
Denote $ \delta_\theta(x_1,\ldots,x_m)=(y_1,\ldots,y_m)$. Then, we have
\[
 x_1=y_1^{-1},
 \quad
 x_i=y_i^{-1}x_{i-1}=(y_1\cdots y_i)^{-1},
 \quad \text{for }2\leq i\leq m.
\]
These expressions define a regular map and its inverse. By the dictionary in \Cref{prop:cocycle-end} and \Cref{thm:inverse-cocycle}, we get \eqref{eq:shift-RB}. When $H$ is connected and reductive, \Cref{thm:main-moduli} identifies the stabilizer of the point
$\theta$ (equivalently, of $B$) with the centralizer of
\[
 \theta(G)=\{e\}\times H^{m-1},
\]
which is exactly $H\times Z(H)^{m-1}$. 
\end{proof}

\subsection{Vector groups}

Let $G=\Ga^r$ with $r\geq1$. Since $G$ is commutative,
a Rota--Baxter operator on $G$ is exactly an algebraic group
endomorphism. In characteristic zero, the additive polynomial
identity forces every such map to be linear, i.e.,
\[
B_A(x)=Ax,\quad \text{for some }A\in M_r(\bfk).
\]
The same argument also works for an arbitrary base scheme. Thus, the full functor of Rota--Baxter operators on $G$ is represented by an affine space of dimension $r^2$, and the subfunctor consisting of invertible ones is represented by $\GL_r$. Every Rota--Baxter operator on the abelian Lie algebra $\mathfrak{g}$ algebraizes uniquely to one of these maps. Since the conjugation action of $G=\Ga^r$ is trivial, the quotient stack and the set of its isomorphism classes over $\bfk$ are, respectively,
\[\cM_{\RB}^{\times}(\Ga^r)\simeq \GL_r\times B\Ga^r,\quad   \pi_0(\cM_{\RB}^{\times}(\Ga^r)(\bfk))\simeq \GL_r(\bfk).
\]

However, the projection
$\GL_r\times B\Ga^r\to\GL_r$
is not a good moduli space. To see this, let $V$ be a two-dimensional $\bfk$-vector space with a  basis $\{e_1,e_2\}$. The projection
$\Ga^r\to\Ga$ onto the first factor yields a two-dimensional
representation
\[
  u\cdot e_1=e_1,\quad
  u\cdot e_2=e_2+u_1e_1
  \quad \text{for }u=(u_1,\ldots,u_r).
\]

In the exact sequence of rational $G$-modules
\[
  0\longrightarrow \bfk e_1\longrightarrow V
  \longrightarrow V/\bfk e_1\longrightarrow0,
\]
$V^{\Ga^r}\to(V/\bfk e_1)^{\Ga^r}$ is not
surjective, since the image is zero but $V/\bfk e_1$ is the one-dimensional trivial module of $G$. Hence, the functor of taking invariants is not exact. Tensoring this sequence with $\cO(\operatorname{GL}_r)$ shows that pushforward along the projection $\operatorname{GL}_r\times B\Ga^r\to\operatorname{GL}_r$ is not exact on the category of quasi-coherent sheaves. Hence, this projection is not a good moduli space.

This example contrasts the unobstructed algebraization on vector groups with the lattice obstruction on tori. In addition, it also illustrates the importance of reductivity for good moduli in \Cref{thm:main-moduli}.

\section{Post-algebraic groups, factorization, and rigidity}\label{sec:applications}

Our results in previous sections have several geometric applications. Recall that the abstract and Lie-theoretic post-group and Yang–Baxter constructions are established in \cite{BG,BGST}. We record here their regular algebraic-group realization as post-algebraic groups and connect it to the differentiated modified classical Yang–Baxter operator and the global factorization. Moreover, our rigidity result makes the invertible part of this data fpqc-locally constant modulo gauge over a connected parameter scheme.

\subsection{Post-algebraic groups and post-Lie algebras}

We define a post-algebraic group as follows, which is the regular algebraic-group form of the postgroup in \cite{BGST}.

\begin{definition}
A \emph{post-algebraic group} is an algebraic group $(G,\cdot)$ together with a morphism $\triangleright:G\times G\to G$ such that
each left translation $L_g^\triangleright$ defined by
$L_g^\triangleright(h)=g\triangleright h$ is an automorphism of the algebraic group $(G,\cdot)$, and
\begin{equation}\label{eq:post-group}
 \bigl(g(g\triangleright h)\bigr)\triangleright \ell
 =g\triangleright(h\triangleright\ell),\quad \forall g,h,\ell\in G.
\end{equation}
\end{definition}
Here, an automorphism means that
$g\triangleright(h\ell)=(g\triangleright h)(g\triangleright\ell)$ and $g\triangleright e=e$ as identities of morphisms. Likewise, \eqref{eq:post-group} is also an identity of morphisms. These
conditions remain valid after every base change, and the inverse of the map $L_g^\triangleright$ is also regular for any $g\in G$. The corresponding definition for abstract groups can be found in \cite[Definition 2.1]{BGST}. Its subadjacent product is
\[
 g\circ h=g(g\triangleright h).
\]
Note that \eqref{eq:post-group} makes this product associative. Moreover, it is a group
law by \cite[Theorem 2.4]{BGST}. The following proposition holds for algebraic groups, including regularity in families.

\begin{proposition}\label{prop:post-group}
Let $B:G\to G$ be a Rota--Baxter operator. Put
\begin{equation}\label{eq:post-from-B}
 g\triangleright_Bh=\Ad_{B(g)}(h),
 \quad
 g\circ_Bh=g\Ad_{B(g)}(h).
\end{equation}
Then $(G,\cdot,\triangleright_B)$ is a post-algebraic group.  The inverse
of $g$ for $\circ_B$ is
\begin{equation}\label{eq:circle-inverse}
 g^{\dagger_B}=\Ad_{B(g)^{-1}}(g^{-1}).
\end{equation}
Moreover, the map
\[
 B:(G,\circ_B)\longrightarrow(G,\cdot)
\]
is a homomorphism of algebraic groups, and the two group laws satisfy the
skew-left brace identity
\begin{equation}\label{eq:brace}
 g\circ_B(h\ell)=(g\circ_Bh)g^{-1}(g\circ_B\ell).
\end{equation}
\end{proposition}

\begin{proof}
Put 
\[\Lambda_B:G\times G\longrightarrow G\times G,\quad (g,h)\longmapsto (g,g\triangleright_B h).\]
The morphism $\Lambda_B$ has the regular inverse
\[
 (g,h)\longmapsto\bigl(g,\Ad_{B(g)^{-1}}(h)\bigr),
\]
and every $L_g^{\triangleright_B}=\Ad_{B(g)}$ preserves the group law.
The Rota--Baxter identity
\begin{equation}\label{eq:B-circle-hom}
 B(g\circ_Bh)=B(g)B(h).
\end{equation}
implies that
\[
 (g\circ_Bh)\triangleright_B\ell
 =\Ad_{B(g)B(h)}(\ell)
 =g\triangleright_B(h\triangleright_B\ell),
\]
which is exactly \eqref{eq:post-group}. Meanwhile, \eqref{eq:circle-inverse} follows from direct computation and \eqref{eq:B-circle-hom} implies that $B$ is a homomorphism. Since conjugation by $B(g)$ preserves the multiplication, we have
\[
 g\circ_B(h\ell)
 =g(g\triangleright_Bh)(g\triangleright_B\ell)
 =(g\circ_Bh)g^{-1}(g\circ_B\ell),
\]
which proves \eqref{eq:brace}. The relation with skew-left braces follows from \cite[Proposition 3.22]{BGST}.
\end{proof}

Call a post-algebraic group structure on $(G,\cdot)$ \emph{strictly inner and lifted} if it is
supplied with a regular map $B:G\to G$ such that
\[
 g\triangleright h=\Ad_{B(g)}(h)
\]
and $B:(G,\circ)\to(G,\cdot)$ is a group homomorphism.  The preceding proof shows that Rota--Baxter operators on $G$ are in bijection with strictly inner lifted post-algebraic-group structures on $G$. The map $B$ is invertible as a Rota--Baxter operator exactly when $B:(G,\circ)\to(G,\cdot)$ is an isomorphism of algebraic groups.

However, it is not true that any post-algebraic group is strictly inner and lifted. For example, on $\Ga^2$, the operation
\begin{equation}\label{eq:noninner-post}
 (x,y)\triangleright(u,v)=(u,v+xu)
\end{equation}
is algebraic and satisfies \eqref{eq:post-group}. Furthermore, its subadjacent law is
\[
 (x,y)\circ(u,v)=(x+u,y+v+xu).
\]
However, it cannot have the form \eqref{eq:post-from-B}, since every inner automorphism of $\Ga^2$ is trivial. On the other hand, the assignment from a Rota--Baxter operator to a post-algebraic-group structure is in general not injective. For a torus, since conjugation is trivial, all Rota--Baxter operators induce the same post-algebraic-group structure even though the invertible parameter set in \eqref{eq:torus-D} can be infinite.

Recall that a \emph{post-Lie algebra} is a Lie algebra $(\mathfrak{g},[\cdot,\cdot])$ equipped with a binary product $\triangleright:\mathfrak g\otimes \mathfrak g\rightarrow \mathfrak g$ such that for any $x,y,z\in \mathfrak g$,
\begin{equation}\label{eq:postlie}
\begin{aligned}
	x\triangleright[y,z]
	&=[x\triangleright y,z]+[y,x\triangleright z],\\
	([x,y]+x\triangleright y-y\triangleright x)\triangleright z
	&=x\triangleright(y\triangleright z)
	-y\triangleright(x\triangleright z).
\end{aligned}
\end{equation}
The following proposition provides an infinitesimal variant of \Cref{prop:post-group}.
\begin{proposition}\label{prop:post-lie}
Let $B:G\to G$ be a Rota--Baxter operator. Let $R=d_eB$. The product
\begin{equation}\label{eq:post-lie-product}
 x\triangleright_Ry=[R(x),y]
\end{equation}
makes $\mathfrak g$ into a post-Lie algebra. Its subadjacent Lie algebra is $(\mathfrak g, [\ ,\ ]_R)$ with
\begin{equation}\label{eq:subadjacent-bracket}
 [x,y]_R:=[x,y]+[R(x),y]+[x,R(y)],\quad \forall x,y\in \mathfrak g
\end{equation}
which is the Lie algebra of $(G,\circ_B)$. Moreover, the map
\[
 R:(\mathfrak g, [\ ,\ ]_R)\longrightarrow
   (\mathfrak g, [\ ,\ ])
\]
is a homomorphism of Lie algebras.
\end{proposition}

\begin{proof}
For the product in \eqref{eq:post-lie-product}, the first identity in \eqref{eq:postlie} is the Jacobi identity. The Rota--Baxter identity implies
\[
 R([x,y]_R)=[R(x),R(y)],
\]
so that the second identity also follows from the Jacobi identity. This also proves that $R$ is a Lie algebra homomorphism. By differentiating \eqref{eq:post-from-B}, we obtain \eqref{eq:subadjacent-bracket}, which identifies the Lie algebra of $(G,\circ_B)$ with $(\mathfrak g, [\ ,\ ]_R)$. We prove the
algebraic version of \cite[Proposition 4.7]{BGST}.
\end{proof}

The post-algebraic-group structure also leads to a regular Yang--Baxter map. Namely,
\begin{equation}\label{eq:YB-map}
 r_B(g,h)=\left(g\triangleright_Bh,
 (g\triangleright_Bh)^{\dagger_B}\circ_Bg\circ_Bh\right)
\end{equation}
is a nondegenerate solution of the braid equation
\[
 (r_B\times\id)(\id\times r_B)(r_B\times\id)
 =(\id\times r_B)(r_B\times\id)(\id\times r_B).
\]
The formula is the specialization of \cite[Corollary 3.14]{BGST} and already appears in \cite[Corollary 6.4]{MZL}. The nondegeneracy follows from \cite[Theorem 3.12 and Proposition 3.13]{BGST}.
It is an algebraic automorphism of $G\times G$. Indeed, if
$r_B(g,h)=(u,v)$, then
\[
 p=u\circ_Bv,
 \quad g=pu^{-1},
 \quad h=g^{\dagger_B}\circ_Bp,
\]
which gives a regular inverse. For an adjoint simple algebraic group, the unique
invertible Rota--Baxter operator is the inversion map, and \eqref{eq:YB-map} becomes
\[
 r_B(g,h)=(g^{-1}hg,g).
\]
For a torus, it is the flip map $(g,h)\mapsto(h,g)$. From this example, we see that the Rota--Baxter moduli is generally finer than the moduli obtained after forgetting the lift and retaining only the Yang--Baxter map.

\subsection{Global factorization and the modified classical Yang--Baxter equation}

\begin{proposition}\label{prop:factorization}
Let $B:G\rightarrow G$ be a Rota--Baxter operator. Define
\[
 B_-(g)=B(g),
 \quad
 B_+(g)=gB(g).
\]
Both morphisms are algebraic group homomorphisms $(G,\circ_B)\to(G,\cdot)$, and
\begin{equation}\label{eq:global-factorization}
 g=B_+(g)B_-(g)^{-1}.
\end{equation}
The paired homomorphism
\[
 \Phi_B=(B_+,B_-):(G,\circ_B)\longrightarrow G\times G
\]
is a closed immersion. Let $H_B$ denote $\Phi_B(G)$ and
let $d:G\times G\to G$ be the morphism $d(x,y)=xy^{-1}$. Then, the restriction
\[
 d|_{H_B}:H_B\xrightarrow{\sim}G
\]
is the inverse of $\Phi_B$. Thus, every $g\in G$ has a unique factorization
$g=g_+g_-^{-1}$ subject to $(g_+,g_-)\in H_B$.
If $R=d_eB$ and $A=2R+\id_{\mathfrak g}$, then
\begin{equation}\label{eq:mcybe}
 [A(x),A(y)]-A\bigl([A(x),y]+[x,A(y)]\bigr)=-[x,y].
\end{equation}
\end{proposition}

\begin{proof}
The claim for $B_-$ is \eqref{eq:B-circle-hom}.  For $B_+$,
\[
 B_+(g\circ_Bh)
 =g\Ad_{B(g)}(h)B(g)B(h)
 =gB(g)hB(h).
\]
The factorization is immediate. The map $\Phi_B$ is a section of $d$. Since $d$ is separated, $\Phi_B$ is a closed immersion. Its image is a closed subgroup because $\Phi_B$ is a group homomorphism. The asserted isomorphism and the uniqueness follow from $d\circ\Phi_B=\id_G$.

To prove \eqref{eq:mcybe}, expand its left-hand side. The two terms are
\[
\begin{aligned}
 [A(x),A(y)]
 &=4[R(x),R(y)]+2[R(x),y]+2[x,R(y)]+[x,y],\\
 A([A(x),y]+[x,A(y)])
 &=4[R(x),R(y)]+2[R(x),y]+2[x,R(y)]+2[x,y],
\end{aligned}
\]
where \eqref{eq:lie-rb} is used in the second equation. Note that the difference is exactly
$-[x,y]$.
\end{proof}

Remark that \eqref{eq:mcybe} is the operator form of the modified classical Yang--Baxter equation. Let $\mathfrak g_R$ denote $\mathfrak g$ with the bracket
\eqref{eq:subadjacent-bracket}. By differentiating $B_-$ and $B_+$, we obtain two homomorphisms of Lie algebras
\[
 R_-=R,
 \quad
 R_+=R+\id:\mathfrak g_R\longrightarrow\mathfrak g.
\]
In the normalization of \cite[(12) and Proposition 7]{STS}, the modified operator is $A=2R+\id$, the Lie bracket on the domain $\mathfrak{g}$ is $[x,y]_A=2[x,y]_R$, and the two homomorphisms are
$A-\id=2R_-$ and $A+\id=2R_+$. Thus, these two formulations agree after this common rescaling. The subgroup $H_B$ and its unique compatible factorization are the algebraic form of \cite[Theorem 3.5]{GLS}.

\begin{remark}
Fix a nondegenerate $\mathfrak g$-invariant bilinear form on $\mathfrak g$, which identifies $\mathfrak g^*$ with $\mathfrak g$. Under the additional hypotheses in \cite[Theorem 11]{STS}, we can obtain a local Lax flow from the factorization. 
\end{remark}

\subsection{An application of the rigidity result}

\begin{corollary}\label{cor:application-rigidity}
Assume that $G$ is connected and reductive. Let $S$ be a connected scheme over $\bfk$. For any family of invertible Rota--Baxter operators over a base $S$, the induced map
\[
 S\longrightarrow M_{\RB}^{\times}(G)
\]
is constant. After an fpqc base change $S'\rightarrow S$, the family is obtained from a fixed operator by the gauge action given in \Cref{thm:gauge-all}. Consequently, every family over a connected base has a constant gauge-equivalence class. Moreover, the associated strictly inner lifted post-algebraic group and the factorization data become gauge equivalent to fixed ones after an fpqc base change. 
\end{corollary}

\begin{proof}
By \Cref{thm:main-moduli}, the good moduli space
$M_{\RB}^{\times}(G)$ is discrete. Note that a morphism from a connected scheme to a discrete scheme is constant. The corresponding family therefore lies in a single component of $D_G$. This component is a $G$-orbit isomorphic to $G/C_i$. If
$f:S\to G/C_i$ is the resulting family, set
$S'=S\times_{G/C_i}G$.  Then $S'\to S$ is smooth, surjective, and
quasi-compact, which implies that it is an fpqc cover. The tautological map $S'\to G$ is the required gauge element. Under the dictionary in \Cref{sec:dic}, it carries a fixed representative to the given family over $S'$.
\end{proof}

If $B'=a\cdot B$, the morphism $T_a^B$ is an  isomorphism of algebraic groups
\[
 T_a^B:(G,\circ_B)\xrightarrow{\sim}(G,\circ_{B'}),
\]
and satisfies
\[
 B'_-\circ T_a^B=B_-,
 \quad
 B'_+\circ T_a^B=\Ad_a\circ B_+.
\]
Hence, a gauge equivalence leads to an isomorphism of the lifted factorization data.

\vskip 5pt
\noindent
\textbf{Acknowledgements}
This study began when the second and third authors visited the Chern Institute of Mathematics in Nankai University and attended the conference ``Rota--Baxter Algebras and Related Topics'' hosted by the Institute in July 2025. They would like to thank the Institute for providing an excellent working environment. The first author is partially supported by NSFC (12271265, W2412041,12671037), the Central Government-Guided Fund Project for Local Science and Technology Development in Tianjin (25ZYZYJC00010), the Fundamental Research Funds for the Central Universities and Nankai Zhide Foundation. The third author is partially supported by the Research Grants Council of the Hong Kong Special Administrative Region, China (Project No. 17305223, 17308324).

\vskip 5pt
\noindent
\textbf{Declaration of interests.}
The authors have no conflict of interest to declare.
\vskip 5pt
\noindent
\textbf{Data availability.}
No data were created or analyzed in this study.
\vskip 5pt
\noindent{\bf AI Declaration.} AI tools (ChatGPT by OpenAI) were used to assist with language editing and improve readability. It did not generate the underlying mathematical ideas or proofs. All mathematical content remains the sole responsibility of the authors.


\begin{thebibliography}{BGST24}

\bibitem[Alp13]{Alper}
J. Alper, \emph{Good moduli spaces for {A}rtin stacks}, Ann. Inst. Fourier
  (Grenoble) \textbf{63} (2013), no.~6, 2349--2402. 

\bibitem[BG22]{BG}
V. G. Bardakov and V. Gubarev, \emph{Rota-{B}axter groups, skew left
  braces, and the {Y}ang-{B}axter equation}, J. Algebra \textbf{596} (2022),
  328--351. 

\bibitem[BGST24]{BGST}
C. Bai, L. Guo, Y. Sheng, and R. Tang, \emph{Post-groups,
  $(${L}ie-$)${B}utcher groups and the {Y}ang-{B}axter equation}, Math. Ann.
  \textbf{388} (2024), no.~3, 3127--3167. 

\bibitem[Bri22]{Brion}
M. Brion, \emph{Homomorphisms of algebraic groups: representability and
  rigidity}, Michigan Math. J. \textbf{72} (2022), 51--76. 

\bibitem[Eti03]{EtingofGeom}
P. Etingof, \emph{Geometric crystals and set-theoretical solutions to the
  quantum {Y}ang-{B}axter equation}, Comm. Algebra \textbf{31} (2003), no.~4,
  1961--1973. 

\bibitem[GLS21]{GLS}
L. Guo, H. Lang, and Y. Sheng, \emph{Integration and geometrization of
  {R}ota-{B}axter {L}ie algebras}, Adv. Math. \textbf{387} (2021), Paper No.
  107834, 34. 

\bibitem[Gon21]{Goncharov}
M. Goncharov, \emph{Rota-{B}axter operators on cocommutative {H}opf
  algebras}, J. Algebra \textbf{582} (2021), 39--56. 

\bibitem[Hit87]{Hitchin}
N. Hitchin, \emph{Stable bundles and integrable systems}, Duke Math. J.
  \textbf{54} (1987), no.~1, 91--114. 

\bibitem[Jac62]{Jacobson}
N. Jacobson, \emph{A note on automorphisms of {L}ie algebras}, Pacific J. Math.
  \textbf{12} (1962), 303--315. 

\bibitem[Jan03]{Jantzen}
J. C. Jantzen, \emph{Representations of algebraic groups}, second ed.,
  Mathematical Surveys and Monographs, vol. 107, American Mathematical Society,
  Providence, RI, 2003. 

\bibitem[KZ02]{KoganZelevinsky}
M. Kogan and A. Zelevinsky, \emph{On symplectic leaves and integrable
  systems in standard complex semisimple {P}oisson-{L}ie groups}, Int. Math.
  Res. Not. (2002), no.~32, 1685--1702. 

\bibitem[MZL26]{MZL}
J. Ma, B. Zhang, and J. Liu, \emph{Affine Rota-Baxter groups and
  affine skew braces}, arXiv preprint arXiv:2606.22954 (2026).

\bibitem[STSi83]{STS}
M. A. Semenov-Tyan-Shansky, \emph{What a classical {$r$}-matrix is},
  Funktsional. Anal. i Prilozhen. \textbf{17} (1983), no.~4, 17--33.

\bibitem[Ves03]{Veselov}
A. P. Veselov, \emph{Yang-{B}axter maps and integrable dynamics}, Phys. Lett. A
  \textbf{314} (2003), no.~3, 214--221. 

\end{thebibliography}
\end{document}